\documentclass[a4paper,11pt]{article}
\usepackage{amssymb,amsmath,amsthm}
\usepackage{texdraw}
\usepackage{epic,eepic}

\newtheorem{thm}{Theorem}[section]
\newtheorem{prop}[thm]{Proposition}

\theoremstyle{remark}
\newtheorem{rem}[thm]{Remark}
\newtheorem{expl}[thm]{Example}

\begin{document}

\title{Yang-Baxter maps and symmetries of integrable equations on quad-graphs}

\author{Vassilios G. Papageorgiou$^{1,2}$, Anastasios G. Tongas$^3$ and Alexander P. Veselov$^{4,5}$ \\ \\
$^1$ {\sl Department of Mathematics, } \\ {\sl National Technical University of Athens, 15780 Zografou, Greece.} \\ 
$^2$ on leave of absence from: \\ {\sl
Department of Mathematics, University of Patras, 26500 Patras, Greece.} \\
{\tt vassilis@math.upatras.gr} \\ \\
$^3$ {\sl Department of Mathematics, University of Patras, 26500 Patras, Greece.} \\ 
{\tt tasos@math.upatras.gr} \\ \\
$^4$ {\sl Loughborough University, Loughborough, Leicestershire LE11 3TU, U.K.} \\ 
$^5$ {\sl Landau Institute for Theoretical Physics, Moscow, Russia.} \\
{\tt A.P.Veselov@lboro.ac.uk} }

\maketitle
\begin{abstract}
A connection between the Yang-Baxter relation for maps and
the multi-dimensional consistency property of integrable equations
on quad-graphs is investigated. The approach is based on the
symmetry analysis of the corresponding equations.
It is shown that the Yang-Baxter variables can be chosen as
invariants of the multi-parameter symmetry groups of the equations.
We use the classification results by Adler, Bobenko and Suris to
demonstrate this method. Some new examples of Yang-Baxter
maps are derived in this way from multi-field integrable equations.
\end{abstract}

\newpage

\section{Introduction} \label{intro}

The quantum Yang-Baxter (YB) equation has its origins in the theory of solvable
models in statistical mechanics \cite{yang,baxter} and the quantum inverse scattering
method \cite{TF}. The fact that this equation has also found many
applications in representation theory, the construction of invariants in knot
theory and that it lies at the foundation of quantum groups, gives to the quantum
YB equation a prominent position among the basic equations in
mathematical physics, see e.g. \cite{lambe}, \cite{fuchs}
and references therein.

In its original form, the quantum YB equation is a relation for a linear
operator $R: V \otimes V \rightarrow V \otimes V$, where $V$ is a vector space.
The relation has the form
\begin{equation}
R^{23}\, R^{13}\, R^{12}\, = R^{12}\, R^{13}\, R^{23}\,, \label{eq:YBrel}
\end{equation}
in $End(V \otimes V \otimes V)$, where $R^{13}$ is meant as the
identity in the second factor of the tensor product $ V \otimes V
\otimes V$ and as $R$ in the first and third factors, and
analogously for $R^{12}$, $R^{23}$. Supposing that $X$ is any set,
the maps $R$ from the Cartesian product $X\times X$ into itself,
which satisfy the relation (\ref{eq:YBrel}) of composite maps are
called set theoretic solutions of the quantum YB relation. The
study of set theoretic solutions of the quantum YB equation was
originally suggested by Drinfeld \cite{drin} (see also earlier work
by Sklyanin \cite{Skl}, where the first interesting example of such 
solution was found) and since then they have attracted the
interest of many researchers.

More recently, a general theory on the set theoretic solutions to
the YB relation was developed in \cite{eting1} and the
notion of transfer maps, which can be considered as the dynamical
analogues of the monodromy and transfer matrices in the theory of
solvable models in statistical mechanics, was introduced in \cite{ves1}.
In many interesting examples of YB maps \footnote{
Adopting the terminology in \cite{ves1}, set
theoretic solutions to the YB equation will be referred
in the following simply as YB maps.},
such as maps
arising from geometric crystals \cite{eting2},
the set $X$ has the structure of an algebraic variety and $R$ is a
birational isomorphism. The case of $\mathbb{CP}^1\times
\mathbb{CP}^1$ has been recently discussed in
\cite{ABS1} in relation with the classification of the
so-called {\em quadrirational maps}.

In this paper we investigate the relation between the YB
property for maps and the {\em multi-dimensional consistency}
condition for equations on quad-graphs, which is now commonly
accepted as a definition of integrability for such equations (see 
\cite{N, BS, ABS2}). Although the link between these
two notions was known before
(see e.g. concluding remarks in \cite{ABS2})
it was never explored systematically.

Our approach is based on the symmetry analysis of integrable
equations on quad-graphs. The main idea is that the YB
variables are suitable invariants of their symmetry groups. A good
example is the {\em discrete potential Korteweg -- de Vries} equation
(dpKdV) \cite{hirota, PNC, CN}
\begin{equation}
\label{dpKdV}
(f_{1,2}-f)(f_1-f_2)- \alpha_1+\alpha_2= 0 \,,
\end{equation}
(see notation on Fig. \ref {fig:BSQ}).
It is clearly invariant under the translation
$f \rightarrow f+{\rm const.}$ The invariants
\begin{equation}
x=f_1-f,\quad y=f_{1,2}-f_1,\quad u = f_{1,2} - f_{2},\quad v = f_2 - f_,
\label{eq:ybvarkdv}
\end{equation}
satisfy the relation
\begin{equation}
x+y=u+v. \label{eq:F51}
\end{equation}
and equation (\ref{dpKdV}) is written in terms of them as
\begin{equation}
(x+y)(x-v)=\alpha_1-\alpha_2. \label{eq:F52}
\end{equation}
This allows to express $u,v$ as functions of $x,y$, which leads to the
following YB map
\begin{equation} \label{eq:adlermap}
u = y +\frac{\alpha_1-\alpha_2}{x+y},\qquad
v = x -\frac{\alpha_1-\alpha_2}{x+y}\,,
\end{equation} known as the Adler map \cite{adler}.
Note that the YB variables $x,y,u,v$, are attached to the edges of the lattice.
The fact that the corresponding map satisfies the YB
property follows directly from the 3D consistency property
of dpKdV (see Fig.~\ref{fig:3DYB}).
This construction works for integrable equations
on quad-graphs with one-parameter symmetry group.

One of the main findings of our paper is that this idea works
for multi-parameter symmetry groups if one considers
the extension of the equation on a multi-dimensional lattice.
In that case the edges are replaced by higher dimensional
faces. We show that in such a way one can derive from the same
discrete potential KdV the following YB map
\begin{equation}
\label{Harrison} u = y Q, \quad v = x Q^{-1}, \quad
Q=\frac{(1-\gamma_2) + (\gamma_2 -\gamma_1)\, x + \gamma_2
(\gamma_1-1)\, x\, y} {(1-\gamma_1) + (\gamma_1 -\gamma_2)\, y +
\gamma_1 (\gamma_2-1)\, x\, y}.
\end{equation}
We will call it {\em Harrison} map since it is closely related to
the superposition formula of the B\"acklund transformation for the
Ernst equation in general relativity introduced by
Harrison \cite{Harrison}. After the change of variables $ x\mapsto
1/x,\,\, v \mapsto 1/v, \,\,  y \mapsto y/\gamma_2,\,\, u\mapsto
u/{\gamma_1} $ it coincides with the $F_I$ quadrirational map in
\cite{ABS1}, which corresponds to the most general case
of two conics. Note that the most degenerate case $F_V$ in the
classification of \cite{ABS1} is simply related to the
Adler map.

The plan of the paper is the following. We start in section
\ref{sec2} with the discussion of 3D consistency property for the
equations on quad-graphs. As the examples we choose three
equations from the classification list in \cite{ABS2}.
By considering the
invariants of their one-parameter symmetry groups we derive all
five types of the quadrirational maps from \cite{ABS1}.
Next, in section {\ref{sec3}}, we show how this symmetry method
can be generalized in the case where the lattice equation admits a
multi-parameter symmetry group. This is demonstrated on the
example of the lattice KdV equation by extending it to a three
dimensional cube and using the invariants of a two-parameter
symmetry group as YB variables. Finally, we show how
the Harrison map can be retrieved from the lattice KdV equation by
exploiting its full three-parameter symmetry group and the
consistency property on a four dimensional cube.

In section \ref{sec4} we show that the same idea works equally
well for multi-field integrable lattice equations and we derive
certain examples of multi-component YB maps from lattice
equations in the Boussinesq family, vector Calapso equation and
its specialization to an integrable discrete version of the
$O(n+2)$ nonlinear $\sigma$--model, introduced recently by Schief
\cite{schief}. The paper concludes with perspectives where we
address some questions for future study.

\section{Multi-dimensional consistency and YB maps}
\label{sec2}
\subsection{Equations on quad-graphs and the 3D consistency property}
\label{par21}
\begin{figure}[h]
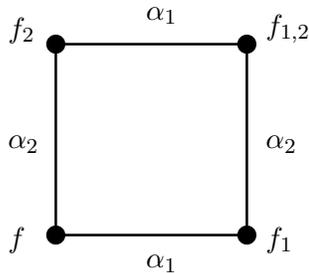

\centertexdraw{\setunitscale 0.5
\linewd 0.03 \arrowheadtype t:F
\move (-1 -1)
\lvec (-1 1) \lvec (1 1) \lvec (1. -1) \lvec(-1 -1)
\move(-1 -1) \fcir f:0.0 r:0.1 \move(-1 1) \fcir f:0.0 r:0.1
\move(1 1) \fcir f:0.0 r:0.1 \move(1 -1) \fcir f:0.0 r:0.1
\htext (-1.5 -1.2) {$f$}
\htext (-1.5 1.) {$f_2$}
\htext (1.2 1) {$f_{1,2}$}
\htext (1.2 -1.2) {$f_{1}$}
\htext (-0.05 1.2) {$\alpha_1$}
\htext (-0.05 -1.4) {$\alpha_1$}
\htext (-1.5 -0.15) {$\alpha_2$}
\htext (1.2 -0.15) {$\alpha_2$}
}
\caption{An elementary quadrilateral.}
\label{fig:BSQ}
\end{figure}

Central to our considerations are integrable discrete equations on
quad-graphs, which are specific equations associated to planar
graphs with elementary quadrilaterals faces. In the simplest case
one has complex fields $f:\mathbb{Z}^2 \rightarrow \mathbb{C}$
assigned on the vertices at sites $(n_1,n_2)$ and two complex lattice
parameters $\alpha_1,\alpha_2$ assigned on the edges of an
elementary square being equal on opposite edges (see Fig.~\ref{fig:BSQ}).
The basic building block of such
equations consists of a relation of the form
\begin{equation}
\mathcal{E}(f,f_{1},f_{2},f_{1,2};\alpha_1,\alpha_2)=0\,, \label{eq:quad}
\end{equation}
between the values of four fields residing on the vertices of each elementary
 quadrilateral for which we use the shorthand notation:
\begin{equation}
f: = f(n_1,n_2), \quad f_{1}:=f(n_1+1,n_2),\quad
f_{2}:=f(n_1,n_2+1), \quad f_{{1,2}}:=f(n_1+1,n_2+1)\,.
\end{equation}

Integrable discrete equations of the above type (\ref{eq:quad})
are listed in a recent classification \cite{ABS2}, where the 3D
consistency property (see below) and some additional conditions
were imposed. From that list we consider the following equations
\begin{eqnarray}
\mathcal{E}_1: & &  (f_{1,2}-f)(f_1-f_2)- \alpha_1+\alpha_2= 0\,,
\label{eq:E1} \\
\mathcal{E}_2: &  &\alpha_1(f\,f_1 + f_2\, f_{1,2} ) - \alpha_2 (f\,f_2 + f_1\, f_{1,2} ) +
\delta ({\alpha_1}^2-{\alpha_2}^2)=0 \,,
\label{eq:E2} \\
\mathcal{E}_3: & &(1-{\alpha_2}^2)(f_1-\alpha_1 f ) (f_2-\alpha_1 f_{1,2}) -
(1-{\alpha_1}^2) (f_2 -\alpha_2 f ) (f_1 - \alpha_2 f_{1,2})=0 \,.
\label{eq:E3}
\end{eqnarray}

Equation $\mathcal{E}_1$ is already mentioned as the
dpKdV equation. Equation $\mathcal{E}_2$ with
$\delta=0$ is the modified discrete KdV or {\em Hirota equation}
\cite{hirota}. If $\delta \neq 0$ we may always assume that
$\delta= 1$ using an appropriate gauge. Equation
$\mathcal{E}_3$ corresponds to the equation labeled as $Q3_{\delta=0}$
in the classification of \cite{ABS2}.
It is contained in the 4-parameter family of the equations derived earlier in 
\cite{frank1}, which contains also discrete versions of potential KdV, 
modified KdV and Schwarzian KdV
(see \cite{frank2} for a more recent discussion).

\begin{figure}[h]
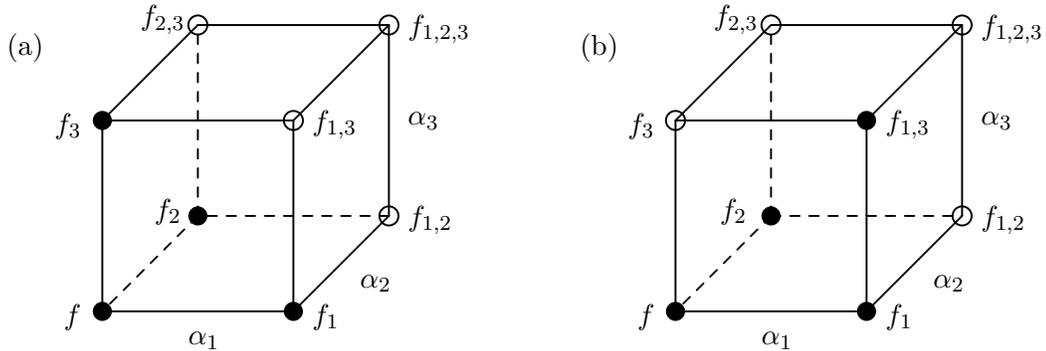

\centertexdraw{ \setunitscale 1. \linewd 0.01 \move (0 0) \linewd
0.01 \lpatt(0.05 0.05) \lvec(0.0 1.0) \lpatt() \lvec (1.0 1.0)
\linewd 0.01 \lvec (1.0 0.0) \lpatt(0.05 0.05) \lvec(0 0)
\lvec(-0.5 -0.5) \lpatt()  \lvec(-0.5 0.5)  \linewd 0.01 \lvec(0
1.0) \move (1. 0)  \linewd 0.01 \lvec(0.5 -0.5) \lvec(-0.5 -0.5)
\move(1.0 1.0) \lvec(0.5 0.5) \lvec(-0.5 0.5) \move(0.5 0.5)
\lvec(0.5 -0.5)

\move(0 0) \fcir f:0.0 r:0.05 \move(-0.5 0.5) \fcir f:0.0 r:0.05
\move(0.5 -0.5) \fcir f:0.0 r:0.05 \move(-0.5 -0.5) \fcir f:0.0
r:0.05 \move(1 1) \lcir r:0.05 \move(0.5 0.5) \lcir r:0.05 \move(0
1) \lcir r:0.05 \move(1 0) \lcir r:0.05

\htext (-1 0.8 ) {(a)} \htext (2 0.8 ) {(b)} \htext (-0.3 0.95)
{$f_{2,3}$} \htext (-0.225 -0.05) {$f_2$} \htext (1.1 0.9)
{$f_{1,2,3}$} \htext (-0.75 0.4) {$f_3$}
\htext (0.6 -0.6) {$f_{1}$} \htext (1.1 -0.1) {$f_{1,2}$} \htext
(-0.7 -0.6) {$f$} \htext (0.6 0.4) {$f_{1,3}$} \htext (0.85 -0.4)
{$\alpha_2$} \htext (-0.05 -0.7) {$\alpha_1$} \htext (1.1 0.45)
{$\alpha_3$}

\move (3 0) \linewd 0.01 \lpatt(0.05 0.05) \lvec (3.0 1.0)  \move
(3 0) \lvec(2.5 -0.5)  \move (3 0) \lvec(4 0) \lpatt() \linewd
0.01 \lvec(4 1.0) \lvec(3 1) \lvec(2.5 0.5)  \lvec(2.5 -0.5)
\lvec(3.5 -0.5) \lvec(4 0) \move(3.5 -0.5) \lvec(3.5 0.5) \lvec(4
1)  \move(3.5 0.5) \lvec(2.5 0.5)

\move(2.5 -0.5) \fcir f:0.0 r:0.05 \move(3.5 -0.5) \fcir f:0.0
r:0.05 \move(4 0) \lcir r:0.05 \move(4 1) \lcir r:0.05 \move(3 0)
\fcir f:0.0 r:0.05 \move(3 1) \lcir r:0.05 \move(2.5 0.5) \lcir
r:0.05 \move(3.5 0.5) \fcir f:0.0 r:0.05

\htext (3.6 -0.6) {$f_{1}$} \htext (4.1 -0.1) {$f_{1,2}$} \htext
(2.3 -0.6) {$f$} \htext (3.6 0.4) {$f_{1,3}$} \htext (3.85 -0.4)
{$\alpha_2$} \htext (2.95 -0.7) {$\alpha_1$} \htext (4.1 0.45)
{$\alpha_3$} \htext (2.7 0.95) {$f_{2,3}$} \htext (2.735 -0.05)
{$f_2$} \htext (4.1 0.9) {$f_{1,2,3}$} \htext (2.25 0.4) {$f_3$}

} \caption{Elementary initial value problems on the cube }
\label{fig:3D}
\end{figure}

The integrability of such equations can be defined using the
{\em three dimensional consistency} property. This means that
the overdetermined system consisting of the difference equations
\begin{eqnarray}
\mathcal{E}(f,f_i,f_j,f_{i,j};\alpha_i,\alpha_j)=0\,, \qquad 1\leq
i < j \leq 3\,, \label{eq:Hij}
\end{eqnarray}
and their shifted versions, is consistent on the three-dimensional
lattice $\mathbb{Z}^3$. In practice, this property is verified as
follows \cite{N,BS, ABS2}. Consider an elementary initial value
problem on the three-dimensional cube with initial data assigned
on four vertices, not all of them lying on the same face. One such
initial configuration is depicted in Fig.~\ref{fig:3D}(a) with
initial values $f,f_i,\,  1 \leq i\leq 3$. Using equations
(\ref{eq:Hij}) on the three faces adjacent to the vertex with
value $f$, we determine uniquely the values $f_{i,j}, \, 1 \leq
i<j \leq 3$, in terms of the initial data. Then using shifted
versions of (\ref{eq:Hij}) on each of the remaining three faces,
we evaluate $f_{1,2,3}$ in three different ways. Consistency means
that one obtains the same value for $f_{1,2,3}$ in terms of the
initial data $f,f_i,\,  1 \leq i\leq 3$ (independent of the way we
choose to evaluate it).
\vspace{0.2cm}
For the dpKdV equation (\ref{dpKdV}) this value  is
\begin{equation}
f_{1,2,3} = \frac{ (\alpha_1-\alpha_2)f_1\,f_2 +
(\alpha_3-\alpha_1) f_1\,f_3 + (\alpha_2-\alpha_3)f_2\,f_3}
{ (\alpha_2-\alpha_1) f_3 + (\alpha_1-\alpha_3) f_2 +
(\alpha_3-\alpha_2)f_1}\,. \label{eq:3dkdv}
\end{equation}
Note that the right hand side of equation  (\ref{eq:3dkdv}) is
invariant under any permutation of the indices $(1,2,3)$ which
label the field variables and the associated lattice parameters.

Another initial data configuration, which is best adapted to the
YB property that we consider, is depicted in Fig.~\ref{fig:3D}(b).
A third possible initial configuration is to give the values
$f,f_1,f_2,f_{1,2,3}$. The latter two configurations are
equivalent to the first one by using the equation on one of the
faces. For example, by using the front face equation we can
exchange the value $f_{1,3}$ to $f_{3}$ in the set of initial data.

Using the fact that dpKdV equation  possesses the 3D consistency property,
one can show in a similar manner that it can be consistently imposed on each
$2$-dimensional face of a $4$-cube. Since we are going to use this property later
on we describe explicitly its derivation.

For given initial values $f,f_i$ $i=1,2,3,4$, we determine the
shifted values of the fields involving any two different
directions, using the equations
\begin{equation}
(f_{i,j}-f) (f_i- f_j) = \alpha_i- \alpha_j\,,
\label{eq:4dkdv}
\end{equation}
$1\leq i < j \leq 4$ (see Fig.~\ref{fig:4D}). Successively, since dpKdV
is $3$-dimensional consistent, we determine the values
$f_{ijk}$,  $1\leq i < j < k\leq 4$. Then the value $f_{1,2,3,4}$ can
be found in six different ways, using the dpKdV equations on the six
2--dimensional facets containing the vertex where the value
$f_{1,2,3,4}$ is assigned. This vertex is contained also in four
cubes, each one of them containing three of the six facets, and
the incidence relations are such that taking into account the
three dimensional consistency on each of the four cubes one proves
that the value $f_{1,2,3,4}$ is uniquely determined in terms of
initial data. By direct calculations also we find that this value
is independent of the way that we used to calculate it, and equals
\begin{equation}
f_{1,i,j,k} = \frac{\underset{ijk}{\boldsymbol{\sigma}}(\alpha_1
\alpha_i f_{i,j} + \alpha_j \alpha_k f_{j,k}) (f_1-f_i)(f_j-f_k)}
{\underset{ijk}{\boldsymbol{\sigma}} (\alpha_1 \alpha_i + \alpha_j
\alpha_k) (f_1-f_i)(f_j-f_k)}\,, \label{eq:4D}
\end{equation}
where $\underset{ijk}{\boldsymbol{\sigma}}$ denotes the cyclic sum
over the subscripts $(i,j,k)=(2,3,4), (4,2,3), (3,4,2)$.
It can be easily checked that $f_{1,2,3,4}$ given by (\ref{eq:4D}),
remains invariant under any permutation of the indices $(1,2,3,4)$,
thus dpKdV is four dimensional consistent.

\begin{center}
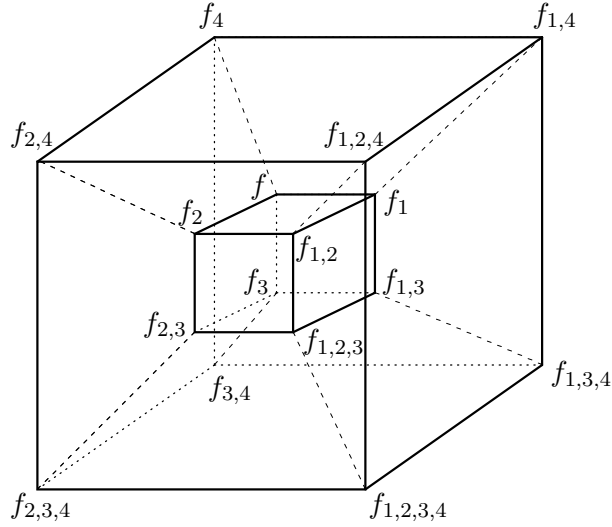
\begin{figure}[h]
\setlength{\unitlength}{0.00037489in}
\begingroup\makeatletter\ifx\SetFigFont\undefined%
\gdef\SetFigFont#1#2#3#4#5{%
  \reset@font\fontsize{#1}{#2pt}%
  \fontfamily{#3}\fontseries{#4}\fontshape{#5}%
  \selectfont}%
\fi\endgroup%
{\renewcommand{\dashlinestretch}{30}
\begin{picture}(7358,6660)(-5000,0)
\path(7110,6390)(7110,1935)
\path(2610,6435)(7110,6435)
\dottedline{100}(180,225)(2610,1935)
\dottedline{100}(3465,4275)(4815,4275)(4815,2925)
    (3465,2925)(3465,4275)
\dottedline{100}(2340,2385)(3465,2925)
\dashline{75.000}(2610,6435)(3465,4275)
\dashline{75.000}(4815,4275)(7110,6435)
\dashline{75.000}(4815,2925)(7110,1935)(7110,1980)
\dashline{75.000}(3690,2385)(4680,225)
\dashline{75.000}(3690,3735)(4680,4725)
\dashline{75.000}(2610,1935)(3465,2925)
\dashline{75.000}(225,4725)(2340,3735)(225,4725)
\dashline{75.000}(180,225)(2340,2385)(180,225)
\thicklines
\path(180,4725)(4680,4725)(4680,225)
    (180,225)(180,4725)
\thinlines
\dottedline{100}(2610,6435)(7110,6435)(7110,1935)
    (2610,1935)(2610,6435)
\thicklines
\path(180,4725)(2610,6435)
\path(2610,6435)(7110,6435)
\path(7110,6435)(7110,1935)
\path(4680,225)(7110,1935)
\path(4680,4725)(7110,6435)
\path(3690,2385)(4815,2925)
\path(4815,4275)(4815,2925)(4815,4275)
\path(3465,4275)(4815,4275)(3465,4275)
\path(2340,3735)(3465,4275)
\path(2340,3735)(3690,3735)(3690,2385)
    (2340,2385)(2340,3735)
\path(3690,3735)(4815,4275)
\put(3100,4275){$f$}
\put(4925,4065){$f_1$}
\put(2385,6635){$f_4$}
\put(3000,2970){$f_3$}
\put(2100,3925){$f_2$}
\put(4905,2970){$f_{1,3}$}
\put(3710,3425){$f_{1,2}$}
\put(-200,5000){$f_{2,4}$}
\put(1625,2385){$f_{2,3}$}
\put(3825,2150){$f_{1,2,3}$}
\put(6975,6635){$f_{1,4}$}
\put(4100,5000){$f_{1,2,4}$}
\put(7220,1710){$f_{1,3,4}$}
\put(4725,-100){$f_{1,2,3,4}$}
\put(2520,1510){$f_{3,4}$}
\put(-200,-100){$f_{2,3,4}$}
\end{picture}
}
\caption{Discrete potential KdV in $\mathbb{Z}^4$}
\label{fig:4D}
\end{figure}
\end{center}

\subsection{YB relation and 3D consistency property}
\label{subsec:YBmaps}

Let $\mathbb{X}$ be any set and $R$ a map of $\mathbb{X} \times
\mathbb{X} $ into itself. Let $R^{ij}:\mathbb{X}^n \rightarrow
\mathbb{X}^n$, where $\mathbb{X}^n=\mathbb{X}\times
\mathbb{X}\times \ldots \times \mathbb{X}$, denotes the map which
acts as $R$ on the $i$ and $j$ factors and as the identity on the
others. More explicitly, let us write $R(x,y)$, $x,y\in
\mathbb{X}$, as
\begin{equation}
R(x,y)=\big(f(x,y),g(x,y)\big).
\end{equation}
Then, for $n \geq 2$ and $1 \leq i,j\leq n$, $i\neq j$ we define
\begin{equation}
R^{ij}(x^1,x^2,\ldots x^n) =
\begin{cases} \,\, (x^1,\ldots,x^{i-1},f(x^i,x^j),x^{i+1}, \ldots,
x^{j-1},g(x^i,x^j),x^{j+1},\ldots x^n) & i<j, \\
\,\, (x^1,\ldots,x^{j-1},g(x^i,x^j),x^{j+1}, \ldots,
x^{i-1},f(x^i,x^j),x^{i+1},\ldots x^n)  & i>j\,.
\end{cases}
\end{equation}
In particular, for $n=2$ we find that $R^{12}=R$ and
$R^{21}(x,y)=\big(g(y,x),f(y,x)\big)$.  The latter
map can be written as a composition of maps as follows
\begin{equation}
R^{21} = P\, R \, P\,,
\end{equation}
where $P$ is the permutation map, i.e. $P(x,y)=(y,x)$.

A map $R$ is called a {\em YB map} if it satisfies the
YB relation (\ref{eq:YBrel}), regarded as an equality of
maps of $\mathbb{X} \times \mathbb{X}\times \mathbb{X}$ into itself.
If in addition the relation $R^{21}\,R=\rm{Id}$ holds, then R is called
reversible YB map.

In a more general setting we may consider a whole family of YB maps
parametrized by continuous parameters $\alpha_i$ rather than a single map.
The YB relation then takes the parameter-dependent form
\begin{equation}
R^{23} (\alpha_2,\alpha_3) \, R^{13} (\alpha_1,\alpha_3) \, R^{12}
(\alpha_1,\alpha_2) = R^{12} (\alpha_1,\alpha_2) \, R^{13}
(\alpha_1,\alpha_3) \, R^{23} (\alpha_2,\alpha_3) \,,
\label{eq:YBcom}
\end{equation}
and the reversibility condition becomes
\begin{equation}
R^{21}(\alpha_2,\alpha_1)\,R(\alpha_1,\alpha_2) = \rm{Id}\,.
\end{equation}

The relation between YB maps and integrable equations on
quad-graphs can be demonstrated in the example of the discrete
potential KdV equation $\mathcal{E}_1$.
As we have already shown in the Introduction, by considering the
differences of the values of the fields assigned on
two adjacent vertices (\ref{eq:ybvarkdv}), we arrive at the Adler map
(\ref{eq:adlermap}).
\begin{figure}[h]
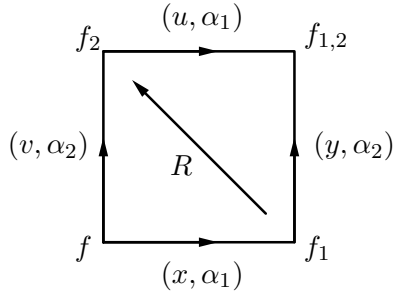

\centertexdraw{\setunitscale 0.5
\linewd 0.03 \arrowheadtype t:F  \arrowheadsize l:0.2 w:0.1
\move (-7 -1)
\lvec (-7 1) \lvec (-5 1) \lvec (-5. -1) \lvec(-7 -1)
\move (-7 -1) \avec( -7 0.1)
\move (-7 -1) \avec( -5.8 -1)
\move (-7 1) \avec( -5.8 1)
\move (-5 -1) \avec( -5 0.1)
\move (-5.3 -0.7) \avec(-6.7 0.7)

\htext (-6.3 -0.3) {$R$}
\htext (-6.4 1.2) {$(u,\alpha_1)$}
\htext (-6.4 -1.5) {$(x,\alpha_1)$}
\htext (-8.0 -0.15) {$(v,\alpha_2)$}
\htext (-4.8 -0.15) {$(y,\alpha_2)$}
\htext (-7.3 -1.2) {$f$}
\htext (-4.9 -1.2) {$f_1$}
\htext (-4.9 1.0) {$f_{1,2}$}
\htext (-7.3 1.0) {$f_2$}

\htext (-3.0 0) {$\phantom{r}$}
\htext (-3.0 1.5) {$\phantom{r}$}
}
\caption{An oriented quadrilateral for the map $R(x,y)=(u,v)$.}
\label{fig:orient}
\end{figure}
There is a different combination for the variables assigned
on the edges of the square, namely
\begin{equation}
x= f \, f_1, \quad y= f_1\, f_{1,2}, \quad u= f_2 f_{1,2}, \quad v= f f_{2} \,.
\label{eq:ybkdv2}
\end{equation}
From the above relations (\ref{eq:ybkdv2}) we deduce that
\begin{equation}
x\, u = y \, v \,.\label{eq:yb2kdv1}
\end{equation}
Moreover, dpKdV can be also written in terms of the variables (\ref{eq:ybkdv2})
as follows
\begin{equation}
y+v-x-u=\alpha_1-\alpha_2 \,. \label{eq:yb2kdv2}
\end{equation}
Solving equations (\ref{eq:yb2kdv1}) and (\ref{eq:yb2kdv2}) for $(u,v)$ we get
the following map
\begin{equation}
u= y \left(1+\frac{\alpha_1-\alpha_2}{x-y}\right),\quad v=
x \left(1+\frac{\alpha_1-\alpha_2}{x-y}\right)
\label{eq:ybmapF4}\,.
\end{equation}

\begin{figure}[h]
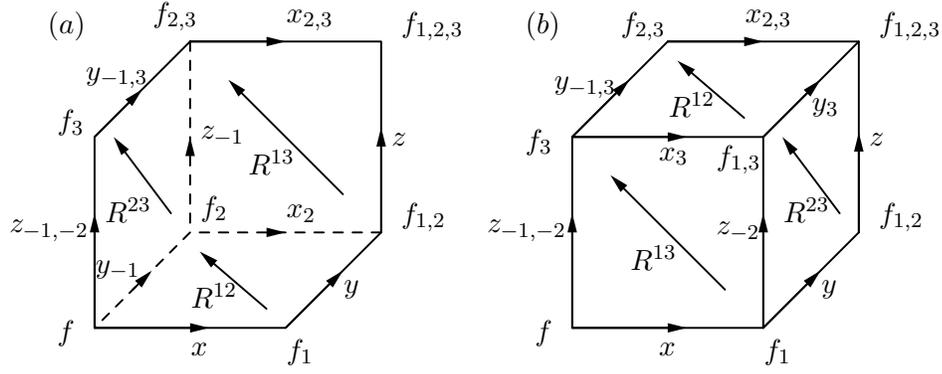

\centertexdraw{ \setunitscale 1. \linewd 0.01 \arrowheadtype t:F
\arrowheadsize l:0.1 w:0.05
\move (1. 0)  \linewd 0.01
\lvec(0.5 -0.5) \lvec(0.5 0.5) \lvec(1.0 1.0)
\move(0.5 0.5) \lvec(-0.5  0.5) \lvec(0.0 1.0)

\move(-2.5 0.0)  \lpatt(0.05 0.05)
\lvec(-3.0 -0.5) \lpatt()
\lvec(-2.0 -0.5)
\lvec(-1.5 0.0)
\lpatt(0.05 0.05) \lvec(-2.5 -0.0)  \lpatt()
\move(-1.5 0.0) \lvec(-1.5 1.0) \lvec(-2.5 1.0)
\lpatt(0.05 0.05) \lvec(-2.5 0.0) \lpatt()
\move(-3.0 -0.5) \lvec(-3.0 0.5) \lvec(-2.5 1.0)
\move(0.5  -0.5) \avec(0.8 -0.2)
\move(0.5  -0.5) \avec(0.5 0.1)
\move(1.0  0.0) \avec(1.0 0.525)
\move(0.5  0.5) \avec(0.8 0.8)
\move(-0.5  0.5) \avec(-0.2 0.8)
\move(0.0  1.0) \avec(0.5 1.0)
\move(-0.5  0.5) \avec(0.1 0.5)
\move(0 1) \lvec(1 1) \lvec(1 0)
\move(-0.5 0.5) \lvec(-0.5 -0.5)  \lvec(0.5 -0.5)
\move(-0.5 -0.5) \avec(-0.5 0.1 )
\move(-0.5 -0.5) \avec(0.1 -0.5 )
\move (-3.0 -0.5) \avec(-2.4 -0.5)
\move (-2.0 -0.5) \avec(-1.7 -0.2)
\move (-2.8 -0.3) \avec(-2.7 -0.2)
\move(-2.1 0.0) \avec(-2.0 0.0)
\move(-2.5 0.4) \avec(-2.5 0.5)
\move(-3.0 -0.5) \avec(-3.0 0.1)
\move(-3.0 0.5) \avec(-2.75 0.75)
\move(-2.5 1.0) \avec(-2.0 1.0)
\move(-1.5 0.0) \avec(-1.5 0.55)
\move(0.9 0.1) \avec(0.6 0.5)
\move(0.3 -0.3) \avec(-0.3 0.3)
\move(0.4 0.6) \avec(0.05 0.9)
\move(-2.6 0.1) \avec(-2.9 0.5)
\move(-1.7 0.2) \avec(-2.3 0.8)
\move(-2.1 -0.4) \avec(-2.45 -0.1)

\htext (-0.25 1.0) {$f_{2,3}$}
\htext (-0.7 -0.6) {$f$}
\htext (1.1 1.0) {$f_{1,2,3}$}
\htext (-0.75 0.4) {$f_3$}
\htext (0.5 -0.7) {$f_{1}$}
\htext (1.1 -0.0) {$f_{1,2}$}
\htext (0.25 0.3) {$f_{1,3}$}
\htext (-0.05 -0.65) {$x$}
\htext (0.8 -0.35) {$y$}
\htext (1.05 0.45) {$z$}
\htext (0.25 -0.05) {$z_{-2}$}
\htext (0.75 0.6) {$y_3$}
\htext (-0.05 0.35) {$x_3$}
\htext (-0.95 -0.05) {$z_{{-1},{-2}}$}
\htext (-0.6 0.7) {$y_{{-1},3}$}
\htext (0.4 1.05) {$x_{2,3}$}

\htext (-3.2 -0.6) {$f$}
\htext (-2.0 -0.7) {$f_1$}
\htext (-1.4 0.0) {$f_{1,2}$}
\htext (-1.4 1.0) {$f_{1,2,3}$}
\htext (-2.7 1.05) {$f_{2,3}$}
\htext (-3.2 0.5) {$f_{3}$}
\htext (-2.45 0.05) {$f_{2}$}
\htext (-2.5 -0.65) {$x$}
\htext (-1.7 -0.35) {$y$}
\htext (-1.45 0.45) {$z$}
\htext (-2.0 0.05) {$x_2$}
\htext (-3.0 -0.25) {$y_{-1}$}
\htext (-2.0 1.05) {$x_{2,3}$}
\htext (-2.45 0.45) {$z_{-1}$}
\htext (-3.05 0.75) {$y_{{-1},3}$}
\htext (-3.45 -0.05) {$z_{{-1},{-2}}$}

\htext (0.6 0.05) {$R^{23}$}
\htext (-0.2 -0.2) {$R^{13}$}
\htext (0.0 0.6) {$R^{12}$}

\htext (-2.5 -0.4) {$R^{12}$}
\htext (-2.2 0.3) {$R^{13}$}
\htext (-2.95 0.05) {$R^{23}$}

\htext (-3.25 1.) {$(a)$}
\htext (-0.75 1.) {$(b)$}

\htext (-0.75 1.3) {$\phantom{R}$}

}
\caption{Three dimensional representation of the YB relation}
\label{fig:3DYB}
\end{figure}

The maps (\ref{eq:adlermap}), (\ref{eq:ybmapF4}), {\em automatically}
satisfy the parameter dependent YB relation (\ref{eq:YBcom}).
Indeed, it is easily shown that the consistency property for a configuration of
initial data on the vertices of a cube, as depicted in Fig.~\ref{fig:3D},
is equivalent to that with initial values $f,f_1,f_{1,2}, f_{1,2,3}$.
These initial data correspond to the values $(x,y,z)$ on the edges
(Fig.~\ref{fig:3DYB}). The 3D consistency property guarantees
that the composite maps
\begin{align}
(a):& \qquad (x,y,z)
\overset{R^{12}}{\longrightarrow}(x_2,y_{-1},z)
\overset{R^{13}}{\longrightarrow}(x_{2,3},y_{-1},z_{-2})
\overset{R^{23}}{\longrightarrow}(x_{2,3},y_{{-1},3},z_{{-1},{-2}})
\\
(b):& \qquad (x,y,z)
\overset{R^{23}}{\longrightarrow}(x,y_3,z_{-2})
\overset{R^{13}}{\longrightarrow}(x_{3},y_3,z_{{-1},{-2}})
\overset{R^{12}}{\longrightarrow}(x_{2,3},y_{{-1},3},z_{{-1},{-2}})
\end{align}
appearing in equation (\ref{eq:YBcom}), applied on $(x,y,z)$
give identical values for $(x_{2,3},y_{{-1},3},z_{{-1},{-2}})$.

Analysing these two examples one notices that the variables
$x,y,u,v$, which we call YB variables, are {\em invariants} of certain
{\em symmetry groups} of the relevant lattice equation.
Now we are going to show that this symmetry method
can be applied in more general situations as well.

\subsection{Lattice invariants of symmetry groups and YB variables}
\label{LieYB}

Let us first recall the basic notions of Lie symmetry methods
applied to lattice equations of the form (\ref{eq:quad}). With
minor modifications these are in accordance with the symmetry
methods applied to algebraic or differential equations
(see e.g. \cite{Olver1} for an extensive study on the subject).

Consider a lattice equation of the form (\ref{eq:quad}) involving one field
$f:\mathbb{Z}^2 \rightarrow \mathbb{C}\, (\mathbb{CP}^1)$.
Let $G$ be a one-parameter group of transformations acting on the domain of the
dependent variables,
\begin{equation}
G: f \mapsto \Phi(n_1,n_2,f;\varepsilon) \,, \quad \varepsilon \in \mathbb{C}\,.
\label{eq:G}
\end{equation}
The prolongation of the group action on the lattice jet space $J$ with coordinates
$(f,f_1,f_2,f_{1,2})$ is specified by
\begin{equation}
G: (f,f_1,f_2,f_{1,2}) \mapsto (\Phi(n_1,n_2,f;\varepsilon),
\Phi_1(n_1+1,n_2,f_1;\varepsilon),
\Phi(n_1,n_2+1,f_2;\varepsilon),
\Phi(n_1+1,n_2+1,f_{1,2};\varepsilon)) \,.\label{eq:Gpr}
\end{equation}
The infinitesimal generator of the group action of $G$
on $f$ is the vector field
\begin{equation}
\mathbf{v} = Q(n_1,n_2,f)\,\partial_f\,, \qquad \mbox{where}\qquad
Q(n_1,n_2,f)=\left.
\frac{\rm d \phantom{\varepsilon}}{{\rm d} \varepsilon}
\Phi(n_1,n_2,f;\varepsilon) \right|_{\varepsilon=0} \,.
\label{eq:vector}
\end{equation}
There is a one-to-one correspondence between connected groups of
transformations and their associated infinitesimal generators since the group
action is reconstructed by the flow of the vector field $\mathbf{v}$ by
exponentiation
\begin{equation}
\Phi(n_1,n_2,f;\varepsilon)=\exp (\varepsilon \, \mathbf{v}) f\,.
\label{eq:orbits}
\end{equation}
The prolongation of the infinitesimal action of $G$ given by (\ref{eq:Gpr}),
is generated by the prolonged vector field
\begin{equation}
\mathbf{\widehat{v}} = Q\,\partial_f +
Q_1\,\partial_{f_1} + Q_2\,\partial_{f_2} +
Q_{1,2}\,\partial_{f_{1,2}}\,,
\label{eq:vectorpr}
\end{equation}
where subscripts denote
$Q_1=Q(n_1+1,n_2,f_1)$, $Q_{1,2}=Q(n_1+1,n_2+1,f_{1,2})$, and so on.

The transformation $G$ is a symmetry of the lattice equation (\ref{eq:quad}),
if it transforms any solution of (\ref{eq:quad}) to another solution of the same
equation. Equivalently, $G$ is a symmetry of equation
(\ref{eq:quad}), if the equation is not affected by the transformation (\ref{eq:Gpr}).
The infinitesimal criterion for $G$ to be a symmetry of equation
(\ref{eq:quad}) is
\begin{equation}
\mathbf{\widehat{v}} \big(\mathcal{E}(f,f_1,f_2,f_{1,2};\alpha_1,\alpha_2)\big) =0\,,
\label{eq:infinv}
\end{equation}
whenever equation (\ref{eq:quad}) holds.

A function $I: J \rightarrow \mathbb{C}$ is a {\em lattice invariant}
of the transformation group $G$, if $I$ is not affected under the action of $G$.
The infinitesimal invariance condition for the lattice invariants is
\begin{equation}
\mathbf{\widehat{v}} (I) = 0 \,. \label{eq:pde}
\end{equation}
Once we have determined a symmetry generator $\mathbf{v}$ of the lattice
equation (\ref{eq:quad}), the corresponding lattice invariants can be
found from the solution of the first order partial differential equation
(\ref{eq:pde}), by using the method of characteristics.
From the corresponding system of ordinary differential equations
we may easily obtain the general solution, since it consists
of equations with separated variables. We assign now to the edges of an
elementary quadrilateral the following YB variables (Fig.~\ref{fig:orient})
\begin{equation}
x=I(f,f_1),\quad y=I(f_1,f_{1,2}),\quad u=I(f_2,f_{1,2}),\quad v=I(f,f_{2})\,,
\label{eq:edgeinv}
\end{equation}
where $I$ is an invariant depending on two neighboring values of $f$.
Since $G$ is a symmetry of the lattice equation, the latter can be written
in terms of these variables:
\begin{equation}
\mathcal{D}(x,y,u,v; \alpha_1,\alpha_2)=0 \,. \label{eq:invlattice}
\end{equation}
This can be done in different ways since the variables (\ref{eq:edgeinv})
are not independent; there exists a relation among them
\begin{equation}
\mathcal{F}(x,y,u,v; \alpha_1,\alpha_2)=0 \,,\label{eq:rel}
\end{equation}
following from the fact that the space of $G$-orbits is three-dimensional.

Solving the system of equations (\ref{eq:invlattice}), (\ref{eq:rel}) for $u,v$
in terms of $x,y$ and assuming that the solution is unique, we obtain a
map $R(x,y)=(u,v)$.
\begin{prop}
If the discrete equation $\mathcal{E}$ satisfies the $3D$ consistency
property, then the map $R(x,y)=(u,v),$ which relates the lattice
invariants (\ref{eq:edgeinv}), satisfies the YB relation.
\end{prop}
The proof follows from Fig.~\ref{fig:3DYB}. Similar considerations
hold for multi-field lattice equations, which give rise to
multi-component YB maps (see section \ref{sec4}).

\begin{expl} Consider the dpKdV equation (\ref{dpKdV}).
Two infinitesimal symmetry generators of the latter equation are
\begin{equation}
\mathbf{v}_1 = \partial_f\,,\qquad
\mathbf{v}_2 = (-1)^{n_1+n_2}\, f\, \partial_f \,,
\end{equation}
(see \cite{TTP}). They generate the symmetry transformations
\begin{equation}
G^1 : f \mapsto f+\varepsilon_1\,,  \qquad
G^2 : f \mapsto f \exp \left( \varepsilon_2 (-1)^{n_1+n_2}\right)\,,  \qquad
\varepsilon_1,\varepsilon_2 \in \mathbb{C}\,,
\end{equation}
respectively. The lattice invariants assigned on the edges of a square for
each one of the above symmetry transformations are the variables
(\ref{eq:ybvarkdv}) and (\ref{eq:ybkdv2}), respectively.
\end{expl}

The consideration of the remaining two equations from the list
(\ref{eq:E1})-(\ref{eq:E3}), with the corresponding symmetry generators
\begin{equation}
\mathbf{v}_{\mathcal{E}_2} = (-1)^{n_1+n_2} f\, \partial_f, \qquad
\mathbf{v}_{\mathcal{E}_3} =  f\, \partial_f\,,
\end{equation}
leads to the results summarized in the following Table \ref{Tab1}.

\begin{table}[h]
\caption{Yang-Baxter maps arising from equations 
$\mathcal{E}_2$ and $\mathcal{E}_3$.} \label{Tab1}
\begin{tabular}{lccc} \hline &&  \\ &
\begin{tabular}{c} Yang-Baxter \\ variables \end{tabular} &
\begin{tabular}{c} Functional relation \\and lattice equation \end{tabular} &
\begin{tabular}{c} Yang-Baxter map \end{tabular} \\ && \\  \hline
&& \\ $\mathcal{E}_2$ &
\begin{tabular}{l} $x=f\, f_1/{\alpha_1}$ \\ $y=f_1\, f_{1,2}/{\alpha_2}$ \\ 
$ u=f_2\, f_{1,2}/{\alpha_1}$ \\ $v =f\, f_2/{\alpha_2}$ \end{tabular}&
\begin{tabular}{c} $\displaystyle{{\gamma_1}\, x \, u = {\gamma_2}\, y \, v}$\\
$\displaystyle{\gamma_1\,(x+u + \delta) = \gamma_2\,(y+v+\delta) }$ \\ 
\text{where} $\,\,\,\gamma_i =  {\alpha_i}^2$
\end{tabular}&
\begin{tabular}{c} $\displaystyle{u = \frac{y}{\gamma_1}\, 
\frac{\gamma_1(x+\delta) -\gamma_2(y+\delta)}{x-y}} $ \\  \\
$ \displaystyle{v = \frac{x}{\gamma_2}\, 
\frac{\gamma_1(x+\delta) - \gamma_2(y+\delta)}{x-y} } $
\end{tabular} \\
&& \\ $\mathcal{E}_3$ &
\begin{tabular}{l} $x=f_1/(\alpha_1\, f)$ \\ $ y=f_{1,2}/{(\alpha_2\, f_1)} $ \\
$ u= f_{1,2}/({\alpha_1 f_2})$ \\  $v=f_{2}/({\alpha_2 f} )$ \end{tabular}&
\begin{tabular}{c} $\displaystyle{ x \, y =  u \, v}$\\
$\displaystyle{ \frac{1-x^{-1}}{1-v^{-1}} = 
\frac{1-{\gamma_1}}{1-{\gamma_2}}\,\frac{1-{\gamma_2}\, y}{1-{\gamma_1} \, u}}$ 
\\ \text{where} $\,\,\,\gamma_i =  {\alpha_i}^2$
\end{tabular}&
\begin{tabular}{c} $\displaystyle{u = y\, Q, \qquad  v = x \, Q^{-1}}$ \\  \\
$ {Q=\frac{(1-\gamma_2) + (\gamma_2 -\gamma_1)\, x + 
\gamma_2 (\gamma_1-1)\, x\, y}
{(1-\gamma_1) + (\gamma_1 -\gamma_2)\, y + \gamma_1 (\gamma_2-1)\, x\, y}} $
\end{tabular} \\
&& \\ \hline 
\end{tabular}
\end{table}

These maps are simply related to quadrilateral maps $F_I - F_{III}$ from the
Adler-Bobenko-Suris list of \cite{ABS1}. Namely, setting $\delta=0$
in the YB map constructed from $\mathcal{E}_2$, we retrieve the
map labeled as $F_{III}$ map in \cite{ABS1}. The case $\delta=-1$
corresponds to the $F_{II}$ map. Finally, the YB map constructed from
$\mathcal{E}_3$ under the transformation
$x\mapsto 1/x,\,\, v \mapsto 1/v, \,\,  y \mapsto y/\gamma_2,\,\,
u\mapsto u/{\gamma_1} $ turns into
\begin{equation}
u = \gamma_1\, y\, \widetilde{Q} , \qquad v = \gamma_2\, x\, \widetilde{Q},
\qquad
\widetilde{Q}=\frac{(1-\gamma_2) x + \gamma_2 -\gamma_1 + (\gamma_1-1) y}
{\gamma_2 (1-\gamma_1) x + (\gamma_1 -\gamma_2) x\,y +\gamma_1 (\gamma_2-1) y}\,,
\label{eq:F1abs}
\end{equation}
which corresponds to the $F_I$ map in the classification in \cite{ABS1}.
As we have already mentioned it is closely related to the superposition
formula of the B\"acklund transformation for the
Ernst equation in general relativity introduced by Harrison
\cite{Harrison}, c.f. \cite{TTX}.

The remaining maps $F_{IV}$ and $F_V$ from \cite{ABS1} are related
in a simple way to the maps (\ref{eq:adlermap}) and (\ref{eq:ybmapF4}) derived
from dpKdV already in section \ref{subsec:YBmaps}.

Thus all 5 types of quadrilateral maps from the Adler-Bobenko-Suris classification 
\cite{ABS1} are equivalent to the YB maps coming from the integrable equations 
on quad-graphs.

\begin{rem}
We should mention that the equivalence of the quadrilateral maps considered
in \cite{ABS1} allows independent change of variables $x,y,u,v$ and
therefore does not respect the YB property, which is preserved in general
only under the diagonal action of the M\"obius group.
In particular, quadrilateral maps in general do not satisfy the YB
relation (contrary to what one might conclude from
\cite{ABS1}). However sometimes two YB maps are related
by non-diagonal action (see the example of Harrison map and $F_V$
above). The question how many such pairs exist needs further
investigation (see the discussion of this in \cite{V2}).
\end{rem}

\section{Multi-parameter symmetry groups and multi-dimensional consistency}
\label{sec3}

The purpose of this section is to show that the symmetry method described
in the preceding section works equally well for multi-parameter symmetry groups.
The idea is to consider the extension of the equation into many dimensions using
the 3D-consistency property and then prescribe the YB variables not to
the edges but, for example, to higher dimensional faces.

We demonstrate how the method works in the example of dpKdV equation, which
is invariant under the three-parameter symmetry group $G$ with
infinitesimal generators
\begin{equation}
\mathbf{v}_1 = \partial_f,\qquad
\mathbf{v}_2 = (-1)^{n_1+n_2}\, f\, \partial_f \,, \qquad
\mathbf{v}_3 = (-1)^{n_1+n_2}\, \partial_f \,.
\end{equation}
Their commutators are
\begin{equation}
\left[ \mathbf{v}_1, \mathbf{v}_2\right] = \mathbf{v}_3\,,\quad
\left[ \mathbf{v}_2, \mathbf{v}_3\right] = -\mathbf{v}_1\,, \quad
\left[ \mathbf{v}_1, \mathbf{v}_3\right] = 0\,,
\end{equation}
from which it is immediately seen that $\{\mathbf{v}_1,
\mathbf{v}_2, \mathbf{v}_3\}$ span a solvable Lie algebra. This
algebra is actually isomorphic to the Lie algebra of the group
$Iso(\mathbb{R}^{1,1})$ of isometries of Minkowski plane, so $G$
can be considered as the connected component of identity of this
group.

\subsection{Consistency of dpKdV around a $3$-cube and the $F_{III}$ map}

Consider the dpKdV equation imposed on each face of an elementary
cube (see Fig.~\ref{fig:3D}) and the abelian subgroup $H$ of the
full symmetry group $G$ generated by $\{\mathbf{v}_1,
\mathbf{v}_3\}$ (translations of the Minkowski plane). These two
symmetries can be extended to the corresponding system on
$\mathbb{Z}^3$, e.g. $\mathbf{v}_3=(-1)^{n_1+n_2+n_3}\partial_f$ .

We consider now the following invariants of the subgroup $H$
\begin{equation}
x=f_1-f_3 \,, \quad y= f_{1,2}-f_{1,3}\,, \quad u=  f_{1,2}-f_{2,3}\,,
\quad v= f_2-f_3\,, \label{eq:invFIII}
\end{equation}
assigned on four faces of the $3$-cube. Using the dpKdV equations
\begin{equation}
f_{1,2} - f = \frac{\alpha_1-\alpha_2}{f_1-f_2}\,,\quad
f_{1,3} - f = \frac{\alpha_1-\alpha_3}{f_1-f_3}\,,\quad
f_{2,3} - f = \frac{\alpha_2-\alpha_3}{f_2-f_3}\,,
\end{equation}
we easily find that the invariants (\ref{eq:invFIII}) are related by
\begin{equation}
u\, v= x \, y\,,\qquad u-\frac{\beta_1}{x} = y-\frac{\beta_2}{v}\,,
\label{eq:FIIIsecond}
\end{equation}
where $\beta_1=\alpha_1-\alpha_3$, $\beta_2=\alpha_2-\alpha_3$.
Solving the system (\ref{eq:FIIIsecond}) for $(u,v)$ in terms of $(x,y)$
we obtain the map
\begin{equation}
u= y P\,, \qquad v=x P^{-1}\,, \qquad P = \frac{\beta_1 + x\,y}{\beta_2 + x\,y}\,,
\label{eq:YBFIIIsecond}
\end{equation}
which satisfies the YB relation as it can be checked by direct calculations.
This fact is also related to the higher dimensional consistency of dpKdV on
$\mathbb{Z}^4$ as it is explained in the following.

Considering the five initial values $f_3,f_{1,3},f_1,f_{1,2},f_{1,2,4} $
on the vertices of the $4$-cube (see Fig.~\ref{fig:4D}) one can find the
values on all other vertices in a unique way using the dpKdV equation on
each $2$-dimensional face because of its $4$-dimensional consistency.
From these five initial values we form the differences
\begin{equation}
x=f_1-f_3 \,, \quad y= f_{1,2}-f_{1,3}\, \quad z=  f_{1,2,4}-f_{1,2,3}\,,
\end{equation}
which are assigned on (the diagonals of) the 2-dimensional faces of the 4-cube
(Fig.~\ref{fig:4D}).
We note that the value $f_{1,2,3}$ can be expressed already in terms of
$f_3,f_{1,3},f_1,f_{1,2}$ through the 3-dimensional consistency of dpKdV
on the ``inner" cube.
Next we apply successively the map $R: (x,y)\mapsto (u,v)$ given by
(\ref{eq:YBFIIIsecond}) on the $(a)$ ``inner", ``front", ``left" and
$(b)$ ``right", ``back", ``outer" 3-dimensional cubes to obtain
the following composite maps
\begin{align}
(a):& \qquad (x,y,z)
\overset{R^{12}}{\longrightarrow}(x_2,y_{-1},z)
\overset{R^{13}}{\longrightarrow}(x_{2,4},y_{-1},z_{-1})
\overset{R^{23}}{\longrightarrow}(x_{2,4},y_{{-1},4},z_{{-1},{-2}})\,,
\\
(b):& \qquad (x,y,z)
\overset{R^{23}}{\longrightarrow}(x,y_4,z_{-2})
\overset{R^{13}}{\longrightarrow}(x_{4},y_4,z_{{-1},{-2}})
\overset{R^{12}}{\longrightarrow}(x_{2,4},y_{{-1},4},z_{{-1},{-2}})\,.
\end{align}
The fact that the two ways of obtaining the values for
$(x_{2,4},y_{{-1},4},z_{{-1},{-2}})$  lead to identical results,
and thus to the YB property of the map (\ref{eq:YBFIIIsecond}), is
guaranteed by the $4$D consistency. One can notice that the
evolution of the Yang Baxter variables  takes place on two
parallel layers of the $\mathbb{Z}^3$ lattice i.e.
$\{(n_1,n_2,0)\, , \, (n_1,n_2,1)\, , n_1,n_2 \in \mathbb{Z}\}$.
This is reflected to the fact that  there are six out of the eight
$3$-dimensional faces (cubes) of the $4$-cube involved in the
compatibility.

A final comment about the map (\ref{eq:YBFIIIsecond}) is that under the
transformation $x\mapsto - x^{-1}$, $v\mapsto - v^{-1}$,
$u \mapsto \beta_1 u$, $y \mapsto \beta_2 y$ it becomes the first YB
map in Table \ref{Tab1} for $\delta=0$.
Thus the $F_{III}$ map in the classification of \cite{ABS1}
is also retrieved from  dpKdV in $\mathbb{Z}^3$, by using the invariants
of the symmetry subgroup $H$.

\subsection{Consistency of dpKdV around a $4$-cube and the Harrison map}

We are going to show that the Harrison map (\ref{Harrison}) appears in a similar
manner as previously using now the invariants of the {\em full} symmetry group $G$
of dpKdV and extending both the equation and its
symmetry group in $\mathbb{Z}^4$.

The symmetry group $G$ can be naturally extended to the corresponding
system in $\mathbb{Z}^4$. Now let us consider the following invariants
of this group:
\begin{equation}
x=\frac{f_{1}-f_{3}}{f_{2}-f_{3}}, \qquad
v=\frac{f_{1}-f_{4}}{f_{2}-f_{4}}, \qquad
y= \left(\frac{f_{1}-f_{4}}{f_{2}-f_{4}}\right)_3,\qquad
u= \left(\frac{f_{1}-f_{3}}{f_{2}-f_{3}}\right)_4\,,
\label{eq:xy}
\end{equation}
where the subscript $i$ means the shift in the $i$-th direction.
The natural place for them to live in are the corresponding 2-dimensional faces.
Next we derive the relations between these variables.

First of all equations (\ref{eq:4dkdv}) and their forward shifts with respect to
the lattice directions $3$ and $4$ imply that the following relations
\begin{subequations}
\begin{align}
(f_{1}-f_{3})(f_{1}-f_{4})_3 &= (f_{1}-f_{4}) (f_{1}-f_{3})_4 \,, \label{eq:fracs1a} \\
(f_{2}-f_{3})(f_{2}-f_{4})_3 &= (f_{2}-f_{4})(f_{2}-f_{3})_4 \,, \label{eq:fracs1b}\\
(f_{1,3}-f)(f_{1,4}-f)_3 &= (f_{1,4}-f)(f_{1,3}-f)_4  \,, \label{eq:fracs2a} \\
(f_{2,3}-f)(f_{2,4}-f)_3 &= (f_{2,4}-f)(f_{2,3}-f)_4 \,, \label{eq:fracs2b}
\end{align}
\end{subequations}
hold on the ``back" and ``left" 3-cubes of the 4-cube depicted in
Fig.~\ref{fig:4D}.
Dividing memberwise equations (\ref{eq:fracs1a}), (\ref{eq:fracs1b}) and
rearranging terms we get
\begin{equation}
\frac{ f_{1}-f_{3}}{f_{2}-f_{3}} \left( \frac{ f_{1}-f_{4}}{f_{2}-f_{4}}\right)_3 =
\frac{ f_{1}-f_{4}}{f_{2}-f_{4}} \left( \frac{ f_{1}-f_{3}}{f_{2}-f_{3}}\right)_4  \,.
\label{eq:frac12}
\end{equation}
In terms of the variables (\ref{eq:xy}) equation (\ref{eq:frac12}) reads
\begin{equation}
x\,y=u\,v\,. \label{eq:eq1}
\end{equation}
On the other hand, we can rewrite equation (\ref{eq:fracs2a}) in the
equivalent form
\begin{equation}
(f_{1,3}-f) \left(1- \frac{(f_{2,3}-f)_4}{(f_{1,3}-f)_4}\right) =
(f_{1,4}-f) \left(1- \frac{(f_{2,4}-f)_3}{(f_{1,4}-f)_3}\right)\,,
\label{eq:q1}
\end{equation}
which by using equations (\ref{eq:4dkdv}) reads
\begin{equation}
\frac{\alpha_1-\alpha_3}{f_1-f_3}
\left(1- \frac{\alpha_2-\alpha_3}{\alpha_1-\alpha_3}\,
\left(\frac{f_{1}-f_3}{f_{2}-f_3}\right)_4\right) =
\frac{\alpha_1-\alpha_4}{f_1-f_4}
\left(1- \frac{\alpha_2-\alpha_4}{\alpha_1-\alpha_4}\,
\left(\frac{f_{1}-f_4}{f_{2}-f_4}\right)_3\right) \,.
\label{eq:q2}
\end{equation}
Multiplying both terms of equation (\ref{eq:q2}) with
$(\alpha_1-\alpha_2)(f_1-f_2)$ and rearranging terms,
the latter takes the form
\begin{align}
\left(1- \frac{\alpha_2-\alpha_4}{\alpha_1-\alpha_4}\right)
\left(1-\frac{f_2-f_3}{f_1-f_3}\right)
\left(1- \frac{\alpha_2-\alpha_3}{\alpha_1-\alpha_3}\,
\left(\frac{f_{1}-f_3}{f_{2}-f_3}\right)_{4}\right) &\nonumber  = \\
\left(1- \frac{\alpha_2-\alpha_3}{\alpha_1-\alpha_3}\right)
\left(1-\frac{f_2-f_4}{f_1-f_4}\right)
\left(1- \frac{\alpha_2-\alpha_4}{\alpha_1-\alpha_4}\,
\left(\frac{f_{1}-f_4}{f_{2}-f_4}\right)_3\right)&\,.
\label{eq:q3}
\end{align}
Finally, recalling the defining relations of the variables $x,y,u,v$,
equation (\ref{eq:q3}) becomes
\begin{equation}
(1-{\gamma_2})(1-x^{-1}) (1-{\gamma_1} \, u) =
(1-{\gamma_1}) (1-v^{-1}) (1-{\gamma_2}\, y)\,,
\label{eq:eq2}
\end{equation}
where
\begin{equation}
\gamma_1 = \frac{\alpha_{2}-\alpha_{3}}{\alpha_{1}-\alpha_{3}} , \qquad
\gamma_2 = \frac{\alpha_{2}-\alpha_{4}}{\alpha_{1}-\alpha_{4}}\,.
\end{equation}
A similar calculation starting with (\ref{eq:fracs2b}) and using
(\ref{eq:eq1}), delivers the same relation (\ref{eq:eq2}).

\begin{prop}
The invariants (\ref{eq:xy}) of the symmetry group $G$ of the
dpKdV equation extended to $\mathbb{Z}^4$ are related by the YB
map, which coincides with the Harrison map (\ref{Harrison}).
\end{prop}

Indeed, comparing (\ref{eq:eq1}) and (\ref{eq:eq2}) with the
relations from which we obtain the last map in Table \ref{Tab1},
we deduce that $x,y$ and $u,v$ are related by the Harrison map.
The fact that this map satisfies the YB relation can be also
derived from the consistency property and geometry of the
$5$-dimensional lattice.

\section{Multi-component YB maps}
\label{sec4}
In this section we show that YB maps can be constructed equally well
from $3D$ consistent multi-field discrete equations, for which no classification
scheme exploiting the multi-dimensional consistency property has been
obtained yet.

\subsection{YB map from the discrete modified Boussinesq system}

The discrete modified Boussinesq (dmBSQ) equations \cite{frank3} involve two
fields $f,g:\mathbb{Z}^2\rightarrow \mathbb{CP}^1$ and are given by the system
\begin{equation}
{f}_{1,2} = g\,
\frac{\alpha_1\, {f}_2 - \alpha_2 \, {f}_1}{\alpha_1\, {g}_1 - \alpha_2 \, {g}_2}\,,
\qquad
{g}_{1,2} = \frac{g}{f}\,
\frac{\alpha_1\, {f}_1 \,{g}_2 - \alpha_2 \,  {f}_2 \, {g}_1}
{\alpha_1 \, {g}_1 - \alpha_2 \, {g}_2}\,.
\label{eq:bsq}
\end{equation}
Its $3D$ consistency is provided by a lengthy but straightforward calculation that
delivers the symmetric values
\begin{equation} f_{1,2,3} = f \,
\frac{\underset{ijk}{\boldsymbol{\sigma}} \, \alpha_i \alpha_j f_k
(\alpha_i g_i - \alpha_j g_j)} {\underset{ij}{\boldsymbol{\sigma}}
\, \alpha_i \alpha_j (\alpha_i f_i g_j - \alpha_j f_j g_i)}\,,
\qquad g_{1,2,3} = g \, \frac{ \underset{ijk}{\boldsymbol{\sigma}}
\, \alpha_i \alpha_j g_k (\alpha_i f_j - \alpha_j f_i) }
{\underset{ij}{\boldsymbol{\sigma}} \, \alpha_i \alpha_j (\alpha_i
f_i g_j - \alpha_j f_j g_i) }\,,
\end{equation}
with respect to any permutation of the indices $(1,2,3)$. Here
the cyclic sum $\underset{ijk}{\boldsymbol{\sigma}}$ is over the subscripts
$(i,j,k)=(1,2,3), (3,1,2), (2,3,1)$,
and similarly the cyclic sum $\underset{ij}{\boldsymbol{\sigma}}$ is over
$(i,j)=(1,2), (2,3), (3,1)$.
The explicit dependence of $f_{1,2,3}$, $g_{1,2,3}$ on the values $f,g$
implies that dmBSQ does not satisfy the so-called {\em tetrahedron} property,
which is an additional assumpion in the
classification scheme in \cite{ABS2} for one-field discrete equations.

Using the symmetry generators
\begin{equation}
\mathbf{v}_1 = f\, \partial_f \,, \qquad \mathbf{v}_2 = g\, \partial_g \,,
\end{equation}
of the dmBSQ equations, we define as YB variables the following joint
lattice invariants
\begin{align}
x^1 = \frac{f_1}{f},\quad y^1 = \frac{f_{1,2}}{f_{1}},\quad
& u^1 = \frac{f_{1,2}}{f_2},\quad v^1 = \frac{f_2}{f},\\
x^2 = \frac{g_1}{g},\quad y^2 = \frac{g_{1,2}}{g_{1}},\quad
&u^2 = \frac{g_{1,2}}{g_2},\quad v^2 = \frac{g_2}{g}\,.
\end{align}
It is immediately seen that the above equations imply that
\begin{equation}
x^1 y^1=u^1 v^1,\qquad  x^2 y^2=u^2 v^2\,. \label{eq:ybsq1}
\end{equation}
Moreover, the lattice equations (\ref{eq:bsq}) can be expressed in terms of the
above invariants as follows
\begin{equation}
u^1 v^1 = \frac{\alpha_1 v^1-\alpha_2 x^1}{\alpha_1 x^2-\alpha_2 v^2},\qquad
u^2 v^2 = \frac{\alpha_1 x^1 v^2-\alpha_2 v^1 x^2}{\alpha_1 x^2-\alpha_2 v^2}.
\label{eq:ybsq2}
\end{equation}
Finally, solving equations (\ref{eq:ybsq1}), (\ref{eq:ybsq2}) for $(u^i,v^i)$ we
obtain the reversible YB map
\begin{subequations} \label{eq:YBBSQ2}
\begin{eqnarray}
& &u^1 =y^1\, A, \qquad v^1= x^1 \, A^{-1} \qquad
A=\frac{{\alpha_1}^2 \,x^1 + {\alpha_2}^2\, x^1 \,x^2 \,y^1 +
{\alpha_1}\,{\alpha_2}\,x^2\, y^2}
{{\alpha_1}\, {\alpha_2}\, x^1 + {\alpha_1}^2 \,x^1\, x^2 \,y^1 +
{\alpha_2}^2\, x^2 \,y^2 }, \\ \nonumber \\
& & u^2 =y^2 \, B, \qquad v^2= x^2 \, B^{-1} \qquad
B=\frac{{\alpha_1}^2 \,x^1 + {\alpha_2}^2 \,x^1 \,x^2\, y^1 +
{\alpha_1}\,{\alpha_2}\,x^2\, y^2}
{{\alpha_2}^2\, x^1 + {\alpha_1} \, {\alpha_2}\, x^1 \,x^2 \,y^1 +
{\alpha_1}^2\, x^2 \,y^2 }\,.
\end{eqnarray}
\end{subequations}
\subsection{YB map from the discrete potential Boussinesq system}
Discrete potential Boussinesq (dpBSQ) equations is the second member in the
so-called lattice Gel'fand-Dikii hierarchy \cite{Pap1}.
The dpBSQ equations, in the form they were studied recently in \cite{TN},
involve three fields $f,g,h:\mathbb{Z}^2\rightarrow \mathbb{CP}^1$,
and they are given by the  following system
\begin{subequations} \label{eq:pBSQ}
\begin{align}
h_1&=f\,f_1 - g\,, \label{eq:pBSQ1}\\ h_2&=f\,f_2 - g\,,
\label{eq:pBSQ2}\\ h&=f\,f_{1,2} - g_{1,2} -
\dfrac{\alpha_1-\alpha_2}{f_1-f_2}\,. \label{eq:pBSQ3}
\end{align}
\end{subequations}
For the purposes of the present discussion, equations (\ref{eq:pBSQ})
exhibit the interesting feature that the joint invariants
of two symmetry generators are enough to construct a YB map.
In connection with this issue we note
that for an elementary Cauchy problem on a staircase, we should impose
initial values $(f,g,h)$, $(f_1,g_1)$, $(f_2,g_2)$, only.
From these data the values $(h_1,h_2)$ and $(f_{1,2},g_{1,2},h_{1,2})$ are
determined uniquely. In particular, equations  (\ref{eq:pBSQ1}),
(\ref{eq:pBSQ2}) imply that
\begin{equation}
f_{1,2}= \dfrac{g_1-g_2}{f_1-f_2}\,, \label{eq:pBSQ4}
\end{equation}
and subsequently the values $h_{1,2}$ and $g_{1,2}$ are determined
from equation (\ref{eq:pBSQ1})
(or equivalently (\ref{eq:pBSQ2})) and (\ref{eq:pBSQ3}), respectively.

Using the infinitesimal invariance criterion (\ref{eq:infinv}) to determine
the symmetries of equations (\ref{eq:pBSQ}), we find that two particular
symmetry generators are given by the following vector fields
\begin{equation}
\mathbf{v}_1 = \partial_f + f\,\partial_g + f\, \partial_h\,, \quad
\mathbf{v}_2 = \partial_g - \partial_h \,.
\label{eq:sympBSQ}
\end{equation}
They generate the symmetry transformations
\begin{align}
G^1&: (f,g,h)\mapsto (f+\varepsilon_1,g+\varepsilon_1 f +
\dfrac{{\varepsilon_1}^2}{2} , h+
\varepsilon_1 f + \dfrac{{\varepsilon_1}^2}{2})\,, \\
G^2&: (f,g,h) \mapsto (f,g+\varepsilon_2,h-\varepsilon_2) \,,
\end{align}
respectively. We define now as YB variables the following invariants
\label{eq:YBpBSQ}
\begin{equation} \begin{array}{ll}
x^1=f_1-f\,,  & y^1=f_{1,2}-f_1\,, \\
x^2=g_1- g - f (f_1-f)\,, &  y^2 = g_{1,2}-g_1 - f_1 (f_{1,2}-f_1)\,, \\
x^3=h_1- h - f (f_1-f)\,, & y^3 = h_{1,2}-g_1 - f_1 (f_{1,2}-f_1)\,,  \\  \\
u^1=f_{1,2} - f_2\,, & v^1=f_{2}-f\,, \\
u^2=g_{1,2} - g_2 - f_2 (f_{1,2}-f_2)\,, &  v^2 = g_{2}-g - f (f_{2}-f)\,, \\
u^3=h_{1,2} - h_2 - f_2 (f_{1,2}-f_2)\,, & v^3 = h_{2}-h - f (f_{2}-f)\,.
\end{array}
\end{equation}
They are functionally related by
\begin{subequations} \label{eq:fundep}
\begin{align}
u^1+v^1&=x^1+y^1\,, \label{eq:fundep1} \\
u^2+v^2&=x^2+y^2+x^1 y^1-u^1 v^1\,, \label{eq:fundep2} \\
u^3+v^3&=x^3+y^3+x^1 y^1-u^1 v^1\,. \label{eq:fundep3}
\end{align}
\end{subequations}
Moreover, the system of equations formed by
(\ref{eq:pBSQ1})-(\ref{eq:pBSQ3}), (\ref{eq:pBSQ2})-(\ref{eq:pBSQ3}) and
(\ref{eq:pBSQ4}) can be written in terms of the above invariants as follows
\begin{subequations} \label{eq:bsqinYB}
\begin{align}
x^1&=-y^1 + \dfrac{x^2}{x^1-v^1} - \dfrac{v^2}{x^1-v^1}\,,
\label{eq:bsqinYB1} \\
x^3 &= x^2+y^2 + x^1 y^1 + \dfrac{\alpha_1-\alpha_2}{x^1-v^1}\,,
\label{eq:bsqinYB2} \\
v^3 &= u^2+v^2 + u^1 v^1 + \dfrac{\alpha_1-\alpha_2}{x^1-v^1}\,,
\label{eq:bsqinYB3}
\end{align}
\end{subequations}
Solving the system of equations (\ref{eq:fundep}), (\ref{eq:bsqinYB})
for $(u^i,v^i)$ we obtain the following YB map
\begin{equation} \begin{array}{ll}
u^1 = y^1 - (\alpha_1-\alpha_2)\Gamma^{-1}\,, &
v^1= x^1 + (\alpha_1-\alpha_2)\Gamma^{-1}\,, \\
u^2= y^2 + (\alpha_1-\alpha_2)\big(\alpha_1-\alpha_2 - 2 y^1 \Gamma\big)\Gamma^{-2}\,, &
v^2 = x^2 + (\alpha_1-\alpha_2)(x^1+y^1)\Gamma^{-1}\,, \\
u^3= y^3 + (\alpha_1-\alpha_2)\big(\alpha_1-\alpha_2 +(x^1-y^1) \Gamma\big)\Gamma^{-2}\,, &
v^3=x^3\,, \end{array}
\end{equation}
where $\Gamma=x^2-x^3+x^1 y^1 + y^2$.

\subsection{YB map from discrete Calapso equation and nonlinear $\sigma$--model.}

In a recent study on discrete isothermic surfaces, Schief \cite{schief}
introduced the following vector generalization of the dpKdV
\begin{equation}
({\boldsymbol{f}}_{1,2} - {\boldsymbol{f}}) =
\frac{\alpha_1-\alpha_2}{|{\boldsymbol{f}}_1-{\boldsymbol{f}}_2|^2}
( {\boldsymbol{f}}_1-{\boldsymbol{f}}_2)\,,
\label{eq:calapso}
\end{equation}
${\boldsymbol{f}}:\mathbb{Z}^2\mapsto \mathbb{C}^n$, under the name
{\em discrete Calapso equation}.
Equation (\ref{eq:calapso}) is three dimensional consistent since for an initial
value configuration ${\boldsymbol{f}},{\boldsymbol{f}}_i$ as in
Fig.~\ref{fig:3D}(a), one finds that the value ${\boldsymbol{f}}_{1,2,3}$
is given by
\begin{equation}
{\boldsymbol{f}}_{1,2,3} =
\frac{\lambda \, |{\boldsymbol{f}}_2-{\boldsymbol{f}}_3|^2
{\boldsymbol{f}}_1 - \mu \,
|{\boldsymbol{f}}_1-{\boldsymbol{f}}_3|^2 {\boldsymbol{f}}_2 +
\nu \, |{\boldsymbol{f}}_1-{\boldsymbol{f}}_2|^2 {\boldsymbol{f}}_3}
{\lambda \, |{\boldsymbol{f}}_2-{\boldsymbol{f}}_3|^2  - \mu \,
|{\boldsymbol{f}}_1-{\boldsymbol{f}}_3|^2 + \nu \,
|{\boldsymbol{f}}_1-{\boldsymbol{f}}_2|^2 }\,,
\end{equation}
where
\begin{equation}
\lambda = (\alpha_1-\alpha_2)(\alpha_1-\alpha_3)\,, \quad
\mu = (\alpha_1-\alpha_2)(\alpha_2-\alpha_3)\,,\quad
\nu  = (\alpha_1-\alpha_3)(\alpha_2-\alpha_3)\,.
\end{equation}
The consistency property is readily checked since ${\boldsymbol{f}}_{1,2,3}$
is invariant under any permutation of the indices $(1,2,3)$ labeling the field
variables and the corresponding lattice parameters.

The aim now is to construct a YB map from the equation
(\ref{eq:calapso}).
Using the translational invariance of equation (\ref{eq:calapso}),
we define the following YB variables
\begin{equation}
{\boldsymbol{x}}={\boldsymbol{f}}_1-{\boldsymbol{f}},\quad
{\boldsymbol{y}}={\boldsymbol{f}}_{1,2}-{\boldsymbol{f}}_1,
\quad {\boldsymbol{u}} =
{\boldsymbol{f}}_{1,2} - {\boldsymbol{f}}_{2},\quad {\boldsymbol{v}} =
{\boldsymbol{f}}_2 - {\boldsymbol{f}},
\label{eq:ybvarcalapso}
\end{equation}
on the edges of a square, which are related by
\begin{equation}
{\boldsymbol{x}}+{\boldsymbol{y}}={\boldsymbol{u}}+{\boldsymbol{v}}.
\label{eq:ybcalapso1}
\end{equation}
On the other hand, equation (\ref{eq:E1}) can be written in terms of
the variables (\ref{eq:ybvarcalapso}) in the form
\begin{equation}
({\boldsymbol{u}}-{\boldsymbol{y}})=
\frac{\alpha_1-\alpha_2}{|{\boldsymbol{x}}+{\boldsymbol{y}}|^2}
({\boldsymbol{x}}+{\boldsymbol{y}}).
\label{eq:ybcalapso2}
\end{equation}
Hence, equations (\ref{eq:ybcalapso1}), (\ref{eq:ybcalapso2}) deliver
the following reversible YB map
\begin{equation}
{\boldsymbol{u}} =
{\boldsymbol{y}} + \frac{\alpha_1-\alpha_2}{|{\boldsymbol{x}}+
{\boldsymbol{y}}|^2}\,({\boldsymbol{x}}+{\boldsymbol{y}}),\qquad
{\boldsymbol{v}} = {\boldsymbol{x}} - \frac{\alpha_1-\alpha_2}
{|{\boldsymbol{x}}+{\boldsymbol{y}}|^2}\,
({\boldsymbol{x}}+{\boldsymbol{y}})\,.
\label{eq:YBcalapso}
\end{equation}
Moreover, in \cite{schief} it was shown that discrete Calapso equation
(\ref{eq:calapso}) can be specialized to an integrable discrete version of
the $O(n+2)$ nonlinear $\sigma$- model. This reduction is accomplished
by imposing the constraint
\begin{equation}
|\boldsymbol{f}|^2=1\,, \label{eq:con}
\end{equation}
on the discrete Calapso equation (\ref{eq:calapso}).
Since the shifted values of $\boldsymbol{f}$ with respect to any lattice
directions should also satisfy constraint (\ref{eq:con}), equation
(\ref{eq:calapso}) is compatible with this constraint whenever
\begin{equation}
2 \boldsymbol{f}\cdot \boldsymbol{f}_2 - 2 \boldsymbol{f} \cdot
\boldsymbol{f}_1=\alpha_1 - \alpha_2 \, .
\end{equation}
This requirement can be satisfied by taking
\begin{equation}
-2 \boldsymbol{f} \cdot \boldsymbol{f}_1 = \alpha_1\, ,
\qquad -2 \boldsymbol{f}\cdot \boldsymbol{f}_2 = \alpha_2\,.
\end{equation}
In terms of the variables (\ref{eq:ybvarcalapso}), the above
constraints translate to
\begin{equation}
|\boldsymbol{x}|^2=2+\alpha_1\,, \qquad |\boldsymbol{y}|^2=2+\alpha_2\,.
\end{equation}
In view of the previous relations the map (\ref{eq:YBcalapso}) obtains the form
\begin{equation}
\boldsymbol{u} =
\boldsymbol{y} + \frac{|\boldsymbol{x}|^2-|\boldsymbol{y}|^2}
{|\boldsymbol{x}+\boldsymbol{y}|^2}\,(\boldsymbol{x}+\boldsymbol{y}),\qquad
\boldsymbol{v} = \boldsymbol{x} - \frac{|\boldsymbol{x}|^2-|\boldsymbol{y}|^2}
{|\boldsymbol{x}+\boldsymbol{y}|^2}\,(\boldsymbol{x}+\boldsymbol{y})\,.
\label{eq:YBsigma}
\end{equation}
By straightforward calculations one finds that
\begin{equation}
|\boldsymbol{u}|^2= |\boldsymbol{x}|^2\,, \qquad
|\boldsymbol{v}|^2=|\boldsymbol{y}|^2\,.
\end{equation}
Using the above identity it is easily established that the map
$R: (\boldsymbol{x},\boldsymbol{y}) \mapsto (\boldsymbol{u},\boldsymbol{v})$
given by (\ref{eq:YBsigma}) is a reversible YB map.
Up to a permutation this map was first considered by Adler \cite{AdlerPolygon}
in the geometric problem about recuttings of the polygons.

\section{Perspectives}
We have shown how the symmetry analysis of integrable equations on
quad-graphs can be used in order to construct YB maps. In
particular, we derived the Harrison map from the consistently
extended discrete potential Korteweg -- de Vries equation to the
$4$-dimensional lattice. The main question now
is how far this example can be generalized. In particular, for a
given multi-parametric symmetry group, is there a general relation
between the structure of the invariants and the geometry of the
objects which the YB variables are assigned to ? What is the dimension of
the lattice in which the discrete equation should be extended to ?
The analysis of other equations from \cite{ABS2} may
clarify these issues.

\subsection*{Acknowledgements}

This work has been started in May 2005 when one of us (APV) visited Patras
within the Socrates Exchange programme between Loughborough University
and University of Patras.

VGP acknowledges partial support from the programme ``C.
Carath\'eodory" of the Research Committee of the University of
Patras and the hospitality of the Department
of  Mathematics  of the National technical University of Athens. 
The work of AGT was supported by the research grant
Pythagoras B-365-015 of the European Social Fund (EPEAEK II). The
work of APV was partially supported by the European network ENIGMA
(contract MRTN-CT-2004-5652) and ESF programme MISGAM.

\end{document}